\documentclass[12pt]{amsart}
\usepackage{amsfonts, amssymb, latexsym, epsfig, amscd}
\usepackage{hyperref,graphicx}

\newtheorem{Definition-Proposition}[Theorem]{Definition-Proposition}
\newtheorem{Main Conjecture}[Theorem]{Main Conjecture}

\theoremstyle{remark}

\theoremstyle{plain}

\newcommand{\cellsize}{11}
\newlength{\cellsz} 
\newsavebox{\cell}
\sbox{\cell}{\begin{picture}(\cellsize,\cellsize)
\put(0,0){\line(1,0){\cellsize}}
\put(0,0){\line(0,1){\cellsize}}
\put(\cellsize,0){\line(0,1){\cellsize}}
\put(0,\cellsize){\line(1,0){\cellsize}}
\end{picture}}
\newcommand\cellify[1]{\def\thearg{#1}\def\nothing{}%
\ifx\thearg\nothing
\vrule width0pt height\cellsz depth0pt\else
\hbox to 0pt{\usebox{\cell} \hss}\fi%
\vbox to \cellsz{
\vss
\hbox to \cellsz{\hss$#1$\hss}
\vss}}
\newcommand\tableau[1]{\vtop{\let\\\cr
\baselineskip -16000pt \lineskiplimit 16000pt \lineskip 0pt
\ialign{&\cellify{##}\cr#1\crcr}}}
\newcommand{\kellsize}{24}
\newlength{\kellsz} 
\newsavebox{\kell}
\sbox{\kell}{\begin{picture}(\kellsize,\kellsize)
\put(0,0){\line(1,0){\kellsize}}
\put(0,0){\line(0,1){\kellsize}}
\put(\kellsize,0){\line(0,1){\kellsize}}
\put(0,\kellsize){\line(1,0){\kellsize}}
\end{picture}}
\newcommand\kellify[1]{\def\thearg{#1}\def\nothing{}%
\ifx\thearg\nothing
\vrule width0pt height\kellsz depth0pt\else
\hbox to 0pt{\usebox{\kell} \hss}\fi%
\vbox to \kellsz{
\vss
\hbox to \kellsz{\hss$#1$\hss}
\vss}}
\newcommand\ktableau[1]{\vtop{\let\\\cr
\baselineskip -16000pt \lineskiplimit 16000pt \lineskip 0pt
\ialign{&\kellify{##}\cr#1\crcr}}}

\newcommand{\sellsize}{36}
\newlength{\sellsz} 
\newsavebox{\sell}
\sbox{\sell}{\begin{picture}(\sellsize,20)
\put(0,0){\line(1,0){\sellsize}}
\put(0,0){\line(0,1){\sellsize}}
\put(\sellsize,0){\line(0,1){\sellsize}}
\put(0,\sellsize){\line(1,0){\sellsize}}
\end{picture}}
\newcommand\sellify[1]{\def\thearg{#1}\def\nothing{}%
\ifx\thearg\nothing
\vrule width0pt height\sellsz depth0pt\else
\hbox to 0pt{\usebox{\sell} \hss}\fi%
\vbox to \sellsz{
\vss
\hbox to \sellsz{\hss$#1$\hss}
\vss}}
\newcommand\stableau[1]{\vtop{\let\\\cr
\baselineskip -16000pt \lineskiplimit 16000pt \lineskip 0pt
\ialign{&\sellify{##}\cr#1\crcr}}}
\begin{document}
\title{Critique of Hirsch's citation index:\\ a 
Combinatorial Fermi problem}

\author{Alexander Yong}
\address{Dept. of Mathematics, University of Illinois at
Urbana-Champaign, Urbana, IL 61801}

\email{ayong@uiuc.edu}

\date{February 17, 2014}
\maketitle

%\section{Introduction}
\pagestyle{plain}

\section{Introduction} 
\subsection{Overview}
In 2005, physicist
J.~E.~Hirsch \cite{Hirsch} proposed the $h$-index to measure the quality of a researcher's output. 
This metric is the largest integer $n$ such 
that the person has $n$ papers with at least $n$ citations each, 
and all other papers have weakly 
less than $n$ citations. Although the original focus of \emph{loc.~cit.}~was on physicists, the $h$-index is now widely popular. For example, {\sf Google Scholar} 
and the {\sf Web of Science} highlight the $h$-index, among other metrics
such as total citation count, in their profile summaries. 

An enticing point made in \emph{loc. cit.} is that the $h$-index 
is an easy and useful supplement to a citation count ($N_{\tt citations}$), since 
the latter metric may be skewed by a small number of highly cited papers or textbooks. In Hirsch's words:
\begin{quote}
``I argue that two individuals with similar $h$s are comparable in terms of their overall scientific impact, even if their total number of papers or their total number of citations is very different. Conversely, comparing two individuals (of the same scientific age) with a similar number of total papers or of total citation count and very different $h$ values, the one with the higher $h$ is likely to be the more accomplished scientist.''
\end{quote}

It seems to us that users might tend to  eyeball differences of $h$s and citation counts among
individuals during their assessments. Instead, one desires a quantitative baseline for what ``comparable'', ``very different'' and ``similar'' actually mean. 
Now, while this would appear to be a matter for statisticians, we show how textbook combinatorics sheds some light on the relationship between the $h$-index and $N_{\tt citations}$. We present a simple model that raises specific concerns about potential misuses of the $h$-index.

To begin, think of the list of a researcher's citations per paper in decreasing order $\lambda=(\lambda_1\geq \lambda_2\geq\cdots \geq \lambda_{N_{\tt papers}})$ as a \emph{partition} of size $N_{\tt citations}$. Graphically, $\lambda$ is identified 
with its \emph{Young diagram}. For example,
$\lambda=(5,3,1,0) \leftrightarrow \tableau{{\bullet }&{\bullet }&{ \ }&{ \ }&{ \ }\\{\bullet }&{\bullet }&{\ }\\{\ }}$.

A combinatorialist will recognize that the $h$-index of $\lambda$ is the side-length of the {\bf Durfee square} (marked using $\bullet$'s above): this is the largest $h\times h$ square that fits
in $\lambda$. This simple observation is nothing new, and appears in both the bibliometric and combinatorial 
literature, see, e.g., \cite{beyond, Analytic}. In particular, since the
Young diagram of
size $N_{\tt citations}$ with maximum $h$-index is roughly a square, 
we see graphically that $0\leq h\leq \lfloor \sqrt{N_{\tt citations}}\rfloor$.

Next, consider the following question:
\begin{center}
Given $N_{\tt citations}$, what is the estimated range of $h$?
\end{center}
Taking only $N_{\tt citations}$ as input hardly seems like sufficient 
information to obtain a meaningful answer. It is exactly for this reason that we call the question a 
\emph{combinatorial Fermi problem}, by analogy with usual \emph{Fermi problems}; see Section~2. 

Since we assume no prior knowledge, consider each citation profile in an unbiased manner.
That is,  each partition of $N_{\tt citations}$ is chosen with equal probability. In fact, there is a beautiful theory 
concerning the asymptotics of these uniform random partitions. 
This was largely developed by A.~Vershik and his collaborators; see, e.g., the survey \cite{Su}.

Actually, we are interested in ``low'' (practical) values of $N_{\tt citations}$ 
where not all asymptotic results are exactly relevant. Instead, we
use generating series and modern desktop computation 
to calculate the probability that a random $\lambda$ has Durfee 
square size $h$. More specifically, we obtain Table~1 below using the  
{\bf Euler-Gauss identity} for partitions:
\begin{equation}
\label{eqn:EulerGauss}
\prod_{i=1}^{\infty} \frac{1}{1-x^i}=1+\sum_{k\geq 1}\frac{x^{k^2}}{\prod_{j=1}^k (1-x^j)^2}.
\end{equation}
The proof of (\ref{eqn:EulerGauss}) via Durfee squares is 
regularly taught to undergraduate combinatorics students; it is recapitulated in Section~3. 
The pedagogical aims of this note
are elaborated upon in both Sections~2 and~3.

\begin{table}[h]
\centering
\begin{tabular}{|c|c|c|c|c|c|c|c|c|c|}
\hline
$N_{\tt citations}$ & $50$ & $100$ & $200$ & $300$ & $400$ & $500$ & $750$ & $1000$ & $1250$\\
\hline
\text{Interval for $h$} & $[2,5]$ & $[3,7]$ & $[5,9]$ & $[7,11]$ & $[8,13]$ & $[9,14]$ & $[11,17]$ & 
$[13,20]$ & $[15,22]$
\\
\hline
\end{tabular}
\end{table}

\begin{table}[h]
\ \ \ \ \ \ \ \ \ \ \ \ \ \ \ \  \ \ \ \ \ \ \ \ \ \ \ \ \ \ \ \ \ \ \ \ \ 
\begin{tabular}{|c|c|c|c|c|c|c|c|c|c|}
\hline
$1500$ & $1750$ & $2000$ & $2500$ & $3000$ & $3500$ & $4000$ & $4500$ & $5000$ & $5500$\\
\hline
$[17,24]$ & $[18,26]$ & $[20,28]$ & $[22,31]$ & $[25,34]$ & $[27,36]$ & $[29,39]$ & $[31,41]$ & $[34,43]$ & $[35,45]$
\\
\hline
\end{tabular}
\end{table}

\begin{table}[h]
\ \ \ \ \ \ \ \ \ \ \ \ \ \ \ \  \ \ \ \ \ \ \ \ \ \ \ \ \ \ \ \ \ \ \ \ \ 
\begin{tabular}{|c|c|c|c|c|c|c|}
\hline
$6000$ & $6500$ & $7000$ & $7500$ & $8000$ & $9000$ & $10000$\\
\hline
$[36,47]$ & $[37,49]$ & $[39,51]$ & $[40,52]$ & $[42,54]$ & $[44,57]$ & $[47,60]$\\
\hline
\end{tabular}
\caption{Confidence intervals for $h$-index}
\end{table}

The asymptotic result we use, due to E.~R.~Canfield-S.~Corteel-C.~D Savage \cite{Can}, gives 
the mode size of the Durfee square when $N_{\tt citations}\to\infty$.
Since their formula is in line with our computations, even for low $N_{\tt citations}$, we reinterpret their work as the
\[\text{rule of thumb for $h$-index: \ \ \ \ \ } h=\frac{\sqrt{6}\log 2}{\pi}\sqrt{N_{\tt citations}}\approx 0.54\sqrt{N_{\tt citations}}.\] 

The focus of this paper is on mathematicians. For the vast majority of those tested, the actual $h$-index computed using {\sf Mathscinet} or {\sf Google scholar} falls into the confidence intervals. 
Moreover, we found that the rule of thumb is fairly accurate for pure mathematicians. For example, Table~2 shows this for
post-$1998$ Fields medalists.\footnote{Citations pre-$2000$ in {\sf Mathscinet} are not complete. {\sf Google scholar} and 
Thompson Reuters' {\sf Web of Science} also have sources of error. We decided that 
{\sf Mathscinet} was our most complete option for analyzing mathematicians. For relatively recent Fields medalists, the effect of lost citations is
reduced.}

In \cite{Hirsch}
it was indicated that the $h$-index has predictive value for winning the Nobel prize. However,
the relation of $h$ index to the Fields medal is, in our opinion, unclear.
A number of the medalists' $h$ values below are shared (or exceeded) by non-contenders of similar academic age, or with those who have the similar citation counts.  Perhaps, this is reflects a cultural difference between
the mathematics and the scientific communities.

In Section~4, we analyze mathematicians in the National Academy of Sciences, where we
show the correlation between the rule of thumb and actual $h$-indices
is $R=0.94$. After removing book citations, $R=0.95$. We also discuss Abel prize winners and associate professors at three research
universities. 

Ultimately, the reader is encouraged to do checks of the estimates themselves. 

\begin{table}[h]
\centering
\begin{tabular}{|c|c|c|c|c|c|}
\hline
Medalist & Award year & $N_{\tt citations}$ & $h$ & Rule of thumb est. & Confidence interval\\
\hline
\hline
T.~Gowers & $1998$ & $1012$ & $15$ & $17.2$ & $[13,20]$\\
\hline
R.~Borcherds & $1998$ & $1062$ & $14$ & $17.6$ & $[14,21]$\\
\hline
C.~McMullen & $1998$ & $1738$ & $25$ & $22.5$ & $[18,26]$\\
\hline
M.~Kontsevich & $1998$ & $2609$ & $23$ & $27.6$ & $[22,32]$\\
\hline
L.~Lafforgue & $2002$ & $133$ & $5$ & $6.2$ & $[4,8]$\\
\hline
V.~Voevodsky  & $2002$ & $1382$ & $20$ & $20.0$ & $[16,23]$\\
\hline
G.~Perelman & $2006$ & $362$ & $8$ & $10.0$ & $[7,12]$\\
\hline
W.~Werner & $2006$ & $1130$ & $19$ & $18.2$ & $[14,21]$\\
\hline
A.~Okounkov & $2006$ & $1677$ & $24$ & $22.1$ & $[18,25]$\\
\hline 
T.~Tao & $2006$ & $6730$ & $40$ & $44.3$ & $[38,51]$\\
\hline
C.~Ng\^{o} & $2010$ & $228$ & $9$ & $8.2$ & $[5,10]$\\
\hline
E.~Lindenstrauss & $2010$ & $490$ & $12$ & $12.0$ & $[9,14]$\\
\hline
S.~Smirnov & $2010$ & $521$ & $12$ & $12.3$ & $[9,15]$\\
\hline
C.~Villani & $2010$ & $2931$ & $25$ & $29.2$ & $[24,33]$ \\
\hline
\end{tabular}
\caption{Fields medalists $1998-2010$}
\end{table}

We discuss three implications/possible applications of our analysis.
\subsection{Comparing $h$'s when $N_{\tt citations}$'s are very different}
It is understood that $h$-index usually grows with $N_{\tt citations}$. 
However, when are citation counts so different that 
comparing $h$'s is uninformative? 
For example, $h_{\text{Tao}}=40$ ($6,730$ citations) while $h_{\text{Okounkov}}=24$ ($1,677$ citations). The model asserts
the probability of $h_{\text{Okounkov}}\geq 32$ is less than $1$ in 
$10$ million. Note the {\sf Math genealogy project}  has fewer than $200,000$ mathematicians. 

These orders of magnitude predict that no mathematician
with $1,677$ citations has an $h$-index of $32$, even though \emph{technically} it can be as high as $40$. Similarly, one predicts the rarity of pure mathematicians with these citations having ``similar''
$h$-indices (in the pedestrian sense). This is relevant when comparing (sub)disciplines with vastly different typical citation counts. We have a theoretical caution about ``eyeballing''.

\subsection{The rule of thumb and the highly cited}
The model suggests the theoretical behavior of the $h$-index for highly-cited scholars. The extent to which these predictions hold true
is informative. This is true not only for individuals, but for entire fields as well.

Actually, Hirsch defined a proportionality constant $a$ by
$N_{\tt citations}=ah^2$
and remarked, ``I find empirically that $a$ ranges between $3$ and $5$.''
This asserts $h$ is between $\sqrt{1/5}\approx 0.45$ and ${\sqrt{1/3}}\approx 0.58$ times
$\sqrt{N_{\tt citations}}$. 

One can begin to try to understand the similarity between Hirsch's empirical upper bound and the rule of thumb. A conjecture of E.~R.~Canfield (private communication, see Section~3) asserts concentration around the mode Durfee square. Thus, \emph{theoretically}, 
one expects the rule of thumb to be nearly correct for $N_{\tt citations}$ large.

Alas, this is empirically not true, even for pure mathematicians. 
However, we observe something related: $0.54{\sqrt{N_{\tt citations}}}$ 
is higher than the actual $h$ for almost every very highly cited ($N_{\tt citations}> 10,000$) scholar in mathematics, physics, computer science and statistics (among others) we considered. 
On the rare occasion this fails, the estimate is only beat by a small percentage ($<5\%$). The drift in the other direction is often quite large ($50\%$ or more is not unusual in certain areas of engineering or biology). 

Near equality occurs among Abel prize winners.
We also considered all prominent physicists highlighted 
in \cite{Hirsch} (except Cohen and Anderson, due to name conflation in {\sf Web of Science}). The guess is always
an upper bound (on average $14-20\%$ too high). Near equality is met by D.~J.~Scalapino 
($25,881$ citations; $1.00$), C.~Vafa ($22,902$ citations; $0.99$), J.~N.~Bahcall ($27,635$ citations; $0.98$); we have given the ratio $\frac{\text{true $h$}}{\text{estimated $h$}}$.

One reason for highly cited people to have lower than expected
$h$-index is that they tend to have highly cited textbooks. Also, famous academics
often run into the ``Matthew effect'' (e.g., gratuitous citations of their most well-known articles or books).

\subsection{Anomalous $h$-indices} 
More generally, our estimates give a way to flag anomalous $h$-indices of active researchers, i.e., those that
are far outside the confidence interval, or, e.g., those for which the rule of thumb is especially 
inaccurate. 

To see what effect book citations has on our estimates, consider the combinatorialist R.~P.~Stanley. 
Since Stanley has $6,510$ citations, we estimate his index as $43.6$.
However, $h_{\text{Stanley}}=35$, a $20\%$ error. Now, $3,237$ of his citations
come from textbooks. Subtracting these, one estimates his $h$-index as $30.9$
while his revised actual $h$-index is $32$, only a $4\%$ error. This kind of
phenomenon was not uncommon; see the appendix.

For another example, consider T.~Tao's {\sf Google scholar} profile. Since he has $30,053$ citations, the rule of thumb predicts his $h$-index is $93.6$. This is far from his actual $h$-index of $65$. Now, his top
five citations (joint with E.~Candes on compressed sensing) are applied. Removing the papers on this topic leaves $13,942$ citations. His new estimate 
is therefore $63.7$ and his revised $h$-index is $61$. 

In many cases we have looked at, once the ``skewing'' feature of the scholar's profile is removed, the remainder of their profile agrees with the rule of thumb.

\subsection{Conclusions and summary}
Whether it be Fields medalists, Abel prize laureates, job, promotion or grant candidates, clearly, the quality of a researcher cannot be fully measured by numerics.  However, in reality, the $h$-index is used, formally or informally, for comparisons. This paper attempts to provide a theoretical and testable framework to quantitatively understand the limits of such evaluations. 
For mathematicians, the accuracy of the rule of thumb suggests that the differences of $h$ index between two mathematicians is strongly influenced by their respective citation counts. 

While discussion of celebrated mathematicians and their statistics makes for fun coffee shop chatter, 
a serious way that $h$-index comes up in faculty meetings concerns relatively junior mathematicians. Consider a scholar $A$ with $100$ citations and $h$ index of $6$ and a scholar $B$ who has $50$ citations and an $h$-index of $4$. Such numbers are not atypical of math assistant professors going 
up for tenure. Our model predicts $h_A$ to be a little bigger than $h_B$. 
Can one really discern what portion of $h_A-h_B$ is a signal of quality?

The problem becomes larger when $A$ and $B$ are in different subject areas. Citations for major works in applied areas tend to have many more citations than in mathematics. In experimental fields, papers may have many coauthors. Since $h$-index does not account for authorship order, this tends to affect our estimates 
for such subjects.

Pure mathematicians have comparatively fewer coauthors, papers and citations. It is not uncommon for, e.g., solutions to longstanding open problems, to have relatively few citations. Thus an explanatory model for pure mathematicians has
basic reasons for being divergent for some other fields. Yet, if this is the case, can the $h$-index really be used universally? This gives us a theoretical reason
to question whether one can make simple comparisons
across fields, even after a rescaling, as has been suggested in \cite{rescale}. 

\section{Combinatorial Fermi Problems}

\subsection{Usual Fermi problems} Fermi problems are so-named after E.~Fermi,
whose ability to obtain good approximate quantitative answers with little data available is legendary. As an illustration, we use the following example \cite{Cooper}:

\begin{center}
How many McDonalds operate in the United States?
\end{center}
There are $10$ McDonalds in Champaign county, which has a population of about $200,000$. \emph{Assume the number of McDonalds
scales with population.} Since the population of the United States is $300$ million, a ``back-of-the-envelope'' calculation estimates the number of McDonalds at  $15,000$. The actual answer, as of $2012$, is $14,157$.

	Using a simplified assumption like the italicized one above is a feature of a Fermi problem. Clearly, the uniform assumption made is not really correct. However, the focus is on good, fast approximations when more careful answers are either too time consuming to determine, or maybe even impossible to carry out. The approximation can then be used to guide further work to determine more
accurate/better justified answers.

Now, although the estimate is rather close to the actual number, when the estimate is not good, the result is even more interesting, as it helps identify a truly faulty assumption. For instance, analogous analysis predicts that the number of Whole Foods in the United States is $0$. 
Apparently, the presence of that company does not scale by population.

Fermi problems/back-of-the-envelope calculations are a standard part of a physics or engineering education. They are of theoretical value in the construction mathematical models, and of ``real world'' value in professions such as management consulting. However, perhaps because the concept is intrinsically non-rigorous, it is not typically part of a (pure) mathematics curriculum. Specifically, this is true for enumerative combinatorics,
even though the subject's purpose is to count the number of certain objects -- which in the author's experience, many students hope has 
non-theoretical applicability.

\subsection{A combinatorial analogue} By analogy we define a {\bf combinatorial Fermi problem}:
\begin{quote}
Fix $\epsilon>0$. Let $S$ be a finite set of combinatorial objects and $\omega:S\to {\mathbb Z}_{\geq 0}$ be a statistic on $S$. Then we estimate the value of $\omega$ on
any element to be the {\bf confidence interval} $[a,b]$ where the \emph{uniform} probability of picking an element of $S$ outside of this range has probability $<\epsilon$.
\end{quote}

By definition, the (ordinary) {\bf generating series} for the {\bf combinatorial problem} $(S,\omega)$ is defined by 
$G_{(S,\omega)}(x)=\sum_{s\in S}x^{\omega(s)}$.
For any $k$, 
$\#\{s\in S: \omega(s)=k\}=[x^k]G_{(S,\omega)}(x)$, 
i.e., the coefficient of $x^k$
in $G_{(S,\omega)}(x)$. Usually, textbook work involves extracting the coefficient
using formulae valid in the ring of formal power series. However, what is often not emphasized in class is that
this coefficient, and $\#S$ itself can be rapidly extracted using a computer algebra system, allowing for a quick determination of the range $[a,b]$. Since the computer does the work, this is our analogue of a ``back-of-the envelope'' calculation.

	For ``reasonable'' values of $\epsilon$ (such as $\epsilon=2\%$), often the range $[a,b]$ is rather tight. In those cases, there may be a theorem of \emph{asymptotic} concentration near a ``typical'' object. However, even if such theorems are known,
this does not solve the finite problem.

The use of the uniform distribution is a quick way to exactly obtain  
estimates that can be compared with empirical data. Ultimately, it invites the user to consider other probability distributions and more sophisticated statistical analysis (just as one should with the McDonald's example), using e.g., Markov Chain Monte Carlo techniques. 

We mention another combinatorial Fermi problem we have considered elsewhere: the count of the number of indigeneous language families in the Americas \cite{Greenberg}. 
That is a situation where essentially there is no way to know with great certainty the true answer.

\section{The Euler-Gauss identity and its application to the $h$-indices}
We apply the perspective of Section~2 to the $h$-index question, where $S={\sf Par}(n)$
and $\omega:S\to {\mathbb Z}_{\geq 0}$ is the size of a partition's Durfee square.  If ${\sf Par}$ is the set of all partitions and $\sigma:{\sf Par}\to {\mathbb Z}_{\geq 0}$
returns the size of a partition, then the generating series for $({\sf Par},\sigma)$ is
$P(x)=\prod_{i=1}^{\infty} \frac{1}{1-x^i}$.
That is, $\#{\sf Par}(n)=[x^N] P(x)$. A sample textbook reference is \cite{Brualdi}.

Recall the {\bf Euler-Gauss identity} (\ref{eqn:EulerGauss})
from the introduction. The well-known combinatorial proof is that every Young diagram $\lambda$
bijectively decomposes into a triple $(D_{\lambda},R_{\lambda},B_{\lambda})$ where $D_{\lambda}$ is a $k\times k$ square, $R_{\lambda}$ is a Young diagram with at most $k$ rows and $B_{\lambda}$ is a partition with at most $k$ columns. That is, $D_{\lambda}$ is the Durfee square, $R_{\lambda}$ be the
shape to the right of the square and $B_{\lambda}$ to be the shape below it. For example:
\[\lambda=\tableau{{\bullet }&{\bullet }&{ \ }&{ \ }&{ \ }\\{\bullet }&{\bullet }&{\ }\\{\ }}\mapsto
\left(\tableau{{\ }&{\ }\\{\ }&{\ }}, \tableau{{\ }&{\ }&{\ }\\{\ }}, \tableau{{\ }}\right)=(D_{\lambda},
R_{\lambda},B_{\lambda}).\]
The generating series for partitions with at most $k$ columns is directly $\prod_{j=1}^k \frac{1}{1-x^j}$. Since conjugation (the ``transpose'') of shape with at most $k$ rows returns a shape with at most $k$ columns, it follows that the generating series for shapes of the first kind is also $\prod_{j=1}^k \frac{1}{1-x^j}$.

From this argument, we see that the generating series for Young diagrams with Durfee square of size $k$ is
$x^{k^2}\prod_{j=1}^k (1-x^j)^{-2}$. We compute for fixed $h, N_{\tt citations}$:
\[{\sf Prob}(\lambda:|\lambda|=N_{\tt citations}, \text{Durfee square of size $k$})=
\frac{[x^{N_{\tt citations}}]x^{k^2}\prod_{j=1}^k (1-x^j)^{-2}}{{\sf Par}(N_{\tt citations})}.
\]

Often textbook analysis ends at the derivation of (\ref{eqn:EulerGauss}). In a classroom,
using a computer to Taylor expand $\sum_{k=a}^b x^{k^2}\prod_{j=1}^k (1-x^j)^{-2}$,
and comparing the coefficients with the known partition numbers allows the instructor to ``physically'' demonstrate the identity to the student.
Varying $a$ and $b$ shows what range of Durfee square sizes are, e.g., $98\%$ likely to occur for partitions of that size.
Interpreted in terms of our $h$-index problem, these same computations are what gives us Table~1.\footnote{Actually, our computation of ${\sf Par}(N_{\tt citations})$ using  generating series became not so easy on a desktop
machine when $N_{\tt citations}$ is a few thousand. Instead, one could use the Hardy-Ramanujun approximation
${\sf Par}(N_{\tt citations})\sim \frac{1}{4N_{\tt citations}\sqrt{3}}e^{\pi\sqrt{\frac{2N_{\tt citations}}{3}}}$.
Even more precisely, one can use {\sf Wolfram Alpha}, which gives the partition numbers for up to a million, which is well beyond our needs.}

As we state in Section~1, the work of \cite{Can} shows the mode Durfree square size is
$\approx 0.54 \sqrt{N_{\tt citations}}$. E.~R.~Canfield's
concentration conjecture states
that for each $\epsilon>0$, 
\begin{equation}
\label{eqn:Canfield}
\lim_{N_{\tt citations}\to\infty} \frac{\text{$\#$ partitions with $(1-\epsilon)\mu<h<(1+\epsilon)\mu$}}{\#{\sf Par}(N_{\tt citations})}\to 1,
\end{equation}
where $\mu=\frac{\sqrt{6}\log 2}{\pi}\sqrt{N_{\tt citations}}$. This is consistent with Table~1. Further discussion
may appear elsewhere. Also, one would like to examine other distributions on Young diagrams, such as the Plancherel measure, which assigns the shape $\lambda$ the probability 
$(f^\lambda)^2/|\lambda|!$ where $f^\lambda$ is the number of \emph{standard Young tableaux} of shape $\lambda$. 

\section{Further comparisons with empirical data}

\subsection{The National Academy of Science}
We compared our rule of thumb against all $119$ mathematicians of the National Academy of Sciences (see the appendix). 
%\footnote{\url{http://en.wikipedia.org/wiki/List_of_members_of_the_National_Academy_of_Sciences_(Mathematics)}}
The correlation
coefficient is $R=0.93$.
After removing
books (as identified in {\sf Mathscinet}), $R=0.95$. 
A serious concern is that many pre-$2000$ citations are not tabulated in {\sf Mathscinet}. 
Nevertheless, in our opinion, the results are still informative. See the comments in Section~4.4.

\begin{figure}[h]
\centering
\includegraphics[width=130mm]{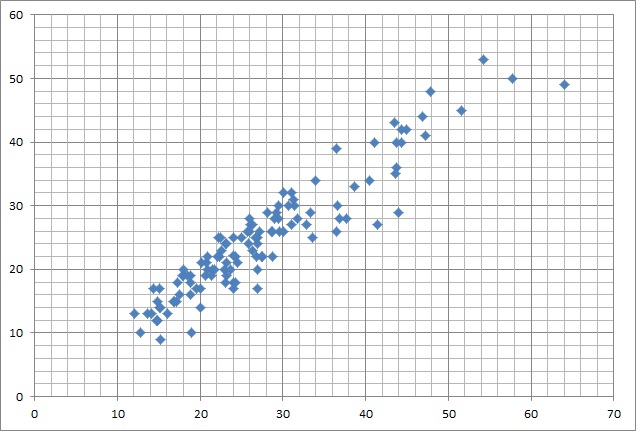}
\caption{Rule of thumb ($x$-axis) versus acutal $h$'s ($y$-axis) for Mathematics members of the National Academy of Sciences}
\end{figure}

\begin{figure}[h]
\centering
\includegraphics[width=130mm]{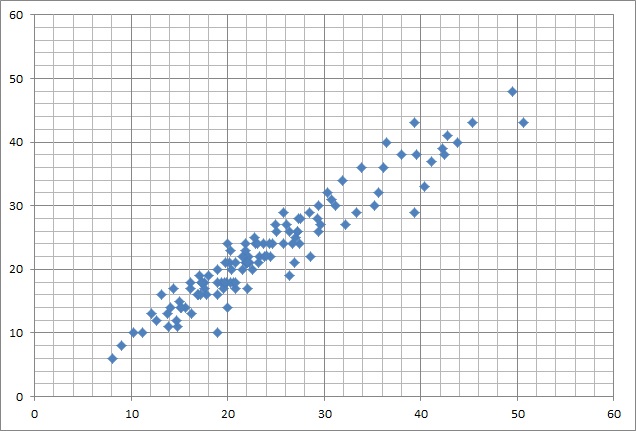}
\caption{Rule of thumb ($x$-axis) versus actual $h$'s ($y$-axis)
for Mathematics members of the National Academy of Sciences (\emph{with books removed})}
\end{figure}

\subsection{Abel prize winners}
Perhaps a closer analogy to the Nobel prize than the Fields medal is the Abel prize, since the latter does not have an age-limit.
The fit with the estimated intervals remains decent; the concern about  pre-$2000$ citations remains.

\begin{table}[h]
\centering
\begin{tabular}{|c|c|c|c|c|c|}
\hline
Laureate & Award year & $N_{\tt citations}$ & $h$ & rule of thumb est. & Estimated range\\
\hline
\hline
J.~P.~Serre & $2003$ & $10119$ & $53$ & $54.3$ & $[47,60]$\\ \hline
I.~Singer & $2004$ & $2982$ & $28$ & $29.5$ & $[24,34]$ \\ \hline
M.~Atiyah & $2004$ & $6564$ & $40$ & $43.7$ & $[37,49]$\\ \hline
P.~Lax & $2005$ & $4601$ & $30$ & $36.6$ & $[31,42]$\\ \hline
L. Carleson & $2006$ & $1980$ & $18$ & $24.0$ & $[19,28]$\\ \hline
S.~R.~S.~Varadhan & $2007$ & $2894$ & $28$ & $29.0$ & $[24,33]$\\ \hline
J.~Thompson & $2008$ & $789$ & $14$ & $15.2$ & $[11,18]$\\ \hline
J. Tits & $2008$ & $3463$ & $28$ & $31.8$ & $[27,36]$\\ \hline
M.~Gromov & $2009$ & $7671$ & $41$ & $47.3$ & $[40,54]$\\ \hline
J.~Tate & $2010$ & $2979$ & $30$ & $29.5$ & $[24,34]$\\ \hline
J.~Milnor & $2011$ & $7856$ & $48$ & $47.9$ & $[41,54]$\\ \hline %38,56
 E.~Szemer\'{e}di  & $2012$  &  $2536$ & $26$ & $27.2$ & $[22,31]$\\ \hline
P.~Deligne & $2013$ & $6567$ & $36$ & $43.8$ & $[37,50]$\\
\hline
\end{tabular}
\caption{Abel prize recipients}
\end{table}
%The $h$-index of highly cited individuals tends to sit in the lower half
%of the estimated interval and the rule of thumb 
%becomes an upper bound. J.~P.~Serre and J.~Milnor are cases where this bound is essentially sharp. %Empirically, for Abel prize winners, the lower bound is around 
%$0.4\sqrt{N_{\tt citations}}$ (for J. Thompson, the least cited Abel prize winner) to $0.47
%\sqrt{N_{\tt citations}}$ (for J.~P.~Serre, our most cited winner).  

\subsection{Associate Professors} Finally, in Table~4 we considered all mathematics associate professors at three research universities. 
Of the $32$ professors, all but five have their $h$-index in the estimated range, and all are at most one
unit outside this range.
 
\begin{table}[h]
\centering
\begin{tabular}{|c|c|c|c|c|}
\hline
& $N_{\tt citations}$ & $h$ & \text{rule of thumb est.} & \text{estimated range}\\
\hline\hline
Department A & \  & \ & \ & \ \\ 
A1 & $19$ & $3$ & $2.4$ & $[1,3]$\\ \hline
A2 & $80$ & $6$ & $4.8$ & $[3,6]$ \\ \hline
A3 & $113$ & $6$ & $5.7$ & $[4,7]$ \\ \hline
A4 & $130$ & $4$ & $6.1$ & $[4,8]$ \\ \hline
A5 & $202$ & $6$ & $7.7$ & $[5,10]$ \\ \hline
A6 & $511$ & $11$ & $12.2$ & $[9,15]$\\ \hline
Department B & \ & \ & \ & \ \\ 
B1 & $30$ & $3$ & $3.0$ & $[1,4]$\\ \hline
B2 & $35$ & $4$ & $3.2$ & $[2,4]$\\ \hline
B3 & $56$ & $4$ & $4.0$ & $[2,5]$\\ \hline
B4 & $56$ & $5$ & $4.0$ & $[2,5]$\\ \hline
B5 & $63$ & $5$ & $4.3$ & $[3,5]$\\ \hline
B6 & $63$ & $6$ & $4.3$ & $[3,5]$\\ \hline
B7 & $78$ & $3$ & $4.8$ & $[3,6]$\\ \hline
B8 & $84$ & $5$ & $4.9$ & $[3,6]$\\ \hline
B9 & $88$ & $7$ & $5.1$ & $[3,6]$\\ \hline
B10 & $122$ & $8$ & $6.0$ & $[4,7]$\\ \hline
B11 & $126$ & $7$ & $6.1$ & $[4,7]$\\ \hline
B12 & $133$ & $6$ & $6.2$ & $[4,8]$\\ \hline
B13 & $133$ & $7$ & $6.2$ & $[4,8]$\\ \hline
B14 & $150$ & $8$ & $6.6$ & $[4,8]$ \\ \hline
B15 & $163$ & $7$ & $6.9$ & $[5,8]$\\ \hline
B16 & $228$ & $10$ & $8.1$ & $[5,10]$\\ \hline
Department C & \ & \ & \ & \ \\
C1 & $10$ & $2$ & $1.7$ & $[1,2]$\\ \hline
C2 & $11$ & $2$ & $1.8$ & $[1,2]$\\ \hline
C3 & $25$ & $3$ & $2.7$ & $[1,3]$ \\ \hline
C4 & $54$ & $4$ & $4.0$ & $[2,5]$ \\ \hline
C5 & $64$ & $5$ & $4.3$ & $[3,5]$ \\ \hline 
C6 & $64$ & $5$ & $4.3$ & $[3,5]$\\ \hline
C7 & $67$ & $6$ & $4.4$ & $[3,5]$ \\ \hline
C8 & $104$ & $6$ & $5.5$ & $[4,7]$ \\ \hline
C9 & $144$ & $8$ & $6.5$ & $[4,8]$ \\ \hline
C10 & $269$ & $5$ & $8.9$ & $[6,11]$ \\
\hline
\end{tabular}
\caption{Associate professors at three research universities}
\end{table}

\subsection{Further study}
It seems to us that the simple model presented describes one force governing $h$-index. However, other forces/sources of noise are at play, depending on
the field or even the fame of the scholar. Future work seeks to better understand this quantitatively, as one works towards more precise models.

The loss of pre-$2000$ citations in {\sf Mathscinet} is significant to how we interpret the results for the National Academy members/Abel prize winners. For example, the rough agreement with the rule of thumb might only reflect an ``equilibrium state'' that arises years after major results have been published. This concern is
partly are allayed by the similar agreement for recent Fields medalists (Table~2). However, as {\sf Mathscinet} reaches further back in tabulating citations, one would try to quantify these effects. In the meantime, use of {\sf Mathscinet} has practical justification since in
promotion and grant decision cases, recent productivity is important. So for these purposes, post-$2000$ data is mostly sufficient.

As a further cross-check, we used the rule of thumb for a broad range of fields using {\sf Google Scholar}.
For scholars with a moderate number of citations, the agreement is often similarly good. Also the rule is an upper bound for the vast majority of
highly cited scholars (but as we have said earlier, much less accurate in some fields). However, these checks have an obvious bias as they only consider people who have set up a profile, so we do not formally
present these results here. 

We propose using the rule of thumb and the confidence intervals as
a basis for a systematic study. We suggest that the rule of thumb reflects an ``ideal scholar''. (This terminology is an allusion to ``ideal gas'' in statistical mechanics.
Indeed, a more conventional use of random partitions concerns the study of Boltzmann statistics on a one-dimensional lattice fermion gas.) 
 Divergence from this ideal is a result of ``anomalies''. For a choice of field, can one statistically distinguish, on quantifiable grounds, scholars who are close to the rule of thumb (in the sense of confidence intervals) from those who are far from it?

\section*{Acknowledgements}
AY thanks George Andrews, Rod Canfield, John D'Angelo, Philippe DiFrancesco, Sergey Fomin, 
Ilya Kapovich, Alexandr Kostochka, Gabriele LaNave, Tom Nevins, Oliver Pechenik, Jim Propp, Bruce Reznick, 
Hal Schenck, Andrew Sills, Armin Straub,
Hugh Thomas, Alexander Woo, Anh Yong and David Yong for helpful comments. We especially thank John D'Angelo and Ilya Kapovich for their
encouragement. This text grew out of combinatorics lectures (Math 413, Math 580) at UIUC; AY thanks the students for their feedback. 
AY was supported by
an NSF grant.

%\appendix
%\section{Data of National Academy of Science Members (Mathematics)}}
% Based on the official list \url{http://nas.nasonline.org/site/Dir/735491561?pg=rslts}
\begin{table}[h]
\caption{(Appendix) Current National Academy of Sciences Members (Mathematics)}
\centering
\begin{tabular}{|c|c|c|c|c|c|c|}
\hline
Member & $N_{\tt citations}$ & Rule of thumb est. & $h$ & non-books only & revised est. & revised $h$\\
\hline
\hline
G. Andrews & $4866$ & $37.7$ & $28$ & $2579$ & $27.4$ & $24$\\ \hline 
M. Artin & $2326$ & $26$ & $26$ & $2097$ & $24.7$ & $24$\\ \hline
M. Aschbacher & $1386$ & $20$ & $17$ & $911$ & $16.3$ & $13$\\ \hline
R. Askey & $2480$ & $26.9$ & $17$ & $1235$ & $19.0$ & $16$\\ \hline
M. Atiyah & $6564$ & $43.7$ & $40$ & $5390$ & $39.6$ & $38$\\ \hline
H. Bass & $2472$ & $26.8$ & $22$ & $1869$ & $23.3$ & $22$\\ \hline
E. Berlekamp & $764$ & $14.9$ & $12$ & $363$ & $10.3$ & $10$\\ \hline
J. Bernstein & $2597$ & $27.5$ & $22$ & $2484$ & $26.9$ & $21$\\ \hline
S. Bloch & $1497$ & $20.9$ & $20$ & $1363$ & $19.9$ & $18$\\ \hline
E. Bombieri & $1746$ & $22.6$ & $23$ & $1608$ & $21.7$ & $22$ \\ \hline
J. Bourgain & $6919$ & $44.9$ & $42$ & $6590$ & $43.8$ & $40$\\ \hline
H. Brezis & $11468$ & $57.8$ & $50$ & $8386$ & $49.5$ & $48$\\ \hline
F. Browder & $2815$ & $28.7$ & $22$ & $2807$ & $28.6$ & $22$\\ \hline
W. Browder & $646$ & $13.7$ & $13$ & $547$ & $12.6$ & $12$\\ \hline
R. Bryant & $1489$ & $20.8$ & $21$ & $1228$ & $18.9$ & $20$\\ \hline
%D. Burkholder & $1089$ & $17.8$ & $17$ & $1089$ & $17.8$ & $17$ \\ \hline
L. Caffarelli & $6745$ & $44.3$ & $42$ & $6280$ & $42.8$ & $41$\\ \hline
E. Calabi & $1224$ & $18.9$ & $18$ & $1224$ & $18.9$ & $18$\\ \hline
L. Carleson & $1980$ & $24$ & $18$ & $1484$ & $20.8$ & $17$\\ \hline
S-Y. Alice Chang & $1828$ & $23.1$ & $24$ & $1806$ & $22.9$ & $24$\\ \hline
J. Cheeger & $3387$ & $31.4$ & $30$ & $3348$ & $31.2$ & $30$\\ \hline
D. Christodoulou & $783$ & $15.1$ & $17$ & $594$ & $13.2$ & $16$ \\ \hline
A. Connes & $6475$ & $43.5$ & $43$ & $5318$ & $39.4$ & $43$\\ \hline
I. Daubechies & $4674$ & $36.9$ & $28$ & $3002$ & $29.6$ & $27$\\ \hline 
P. Deift & $3004$ & $29.6$ & $26$ & $2545$ & $27.2$ & $26$\\ \hline
P. Deligne & $6567$ & $43.8$ & $36$ & $5592$ & $40.4$ & $33$\\ \hline
P. Diaconis & $3233$ & $30.7$ & $30$ & $2970$ & $29.4$ & $30$\\ \hline
S. Donaldson & $2712$ &  $28.1$ &  $29$ & $2277$ & $25.8$ & $29$\\ \hline
E. Dynkin & $1583$ & $21.5$ & $20$ & $1090$ & $17.8$ & $16$\\ \hline
Y. Eliashberg & $1628$ & $21.8$ & $20$ & $1460$ & $20.6$ & $18$\\ \hline
L. Faddeev & $1820$ & $23$ & $20$ & $1285$ & $19.4$ & $18$\\ \hline
C. Fefferman & $3828$ & $33.4$ & $29$ & $3815$ & $33.4$ & $29$\\ \hline
M. Freedman & $1207$ & $18.8$ & $16$ & $990$ & $17$ & $16$\\ \hline
W. Fulton & $5890$ & $41.4$ & $27$ & $1424$ & $20.4$ & $20$\\ \hline
H. Furstenberg & $2064$ & $24.5$ & $21$ & $1650$ & $21.9$ & $21$\\ \hline
D. Gabai & $1314$ & $19.6$ & $17$ & $1314$ & $19.6$ & $17$\\ \hline
%F. Gehring & $1498$ & $20.9$ & $20$ & $1466$ & $20.7$ & $20$\\ \hline
J. Glimm & $1826$ & $23.1$ & $18$ & $1419$ & $20.3$ & $18$\\ \hline
R. Graham & $3881$ & $33.6$ & $25$ & $2280$ & $25.8$ & $24$ \\ \hline
U. Grenander & $895$ & $16.1$ & $13$ & $227$ & $8.1$ & $6$\\ \hline
P. Griffiths & $4581$ & $36.5$ & $26$ & $1692$ & $22.2$ & $22$\\ \hline 
M. Gromov & $7671$ & $47.3$ & $41$ & $6200$ & $42.5$ & $38$\\ \hline
B. Gross & $1692$ & $22.2$  & $25$ & $1635$ & $21.8$ & $24$\\ \hline
\end{tabular}
\end{table}

\begin{table}[b]
\centering
\begin{tabular}{|c|c|c|c|c|c|c|}
\hline
Member & $N_{\tt citations}$ & Rule of thumb est. & $h$ & non-books only & revised est. & revised $h$\\
\hline
\hline
V. Guillemin & $3710$ & $32.9$ & $27$ & $2035$ & $24.4$ & $22$\\ \hline
R. Hamilton & $2490$ & $26.9$ & $20$ & $2392$ & $26.4$ & $19$\\ \hline
%F. Hirzebruch & $1773$ & $22.7$ & $24$ & $1015$ & $17.2$ & $18$\\ \hline
M. Hochster & $1727$ & $22.4$ & $22$ & $1657$ & $22$ & $21$\\ \hline
H. Hofer & $2140$ & $25$ & $25$ & $1928$ & $23.7$ & $24$\\ \hline
MJ. Hopkins & $714$ & $14.4$ & $17$ & $714$ & $14.4$ & $17$\\ \hline
R. Howe & $1680$ & $22.1$ & $22$ & $1579$ & $21.5$ & $22$\\ \hline
H. Iwaniec & $2822$ & $28.7$ & $26$ & $1825$ & $23.1$ & $24$\\ \hline
A. Jaffe & $794$ & $15.2$ & $9$ & $277$ & $9$ & $8$\\ \hline
P. Jones & $1112$ & $18$ & $19$ & $1112$ & $18$ & $19$ \\ \hline
V. Jones & $2025$ & $24.3$ & $18$ & $1669$ & $22.1$ & $17$\\ \hline
R. Kadison & $1922$ & $23.7$ & $20$ & $1042$ & $17.4$ & $18$\\ \hline
R. Kalman & $558$ & $12.8$ & $10$ & $428$ & $11.2$ & $10$\\ \hline
N. Katz & $2370$ & $26.3$ & $23$ & $1582$ & $21.5$ & $20$\\ \hline
D. Kazhdan & $2332$ & $26.1$ & $27$ & $2332$ & $26.1$ & $27$\\ \hline
R. Kirby & $963$ & $16.8$ & $15$ & $678$ & $14.1$ & $14$\\ \hline
S. Klainerman & $2324$ & $26$ & $28$ & $2144$ & $25.0$ & $27$\\ \hline
J. Kohn & $1231$ & $18.9$ & $19$ & $1068$ & $17.6$ & $18$\\ \hline
J. Koll\'{a}r & $3100$ & $30.1$ & $26$ & $1947$ & $23.8$ & $22$\\ \hline
B. Kostant & $2509$ & $27$ & $25$ & $2509$ & $27$ & $25$\\ \hline
R. Langlands & $1466$ & $20.6$ & $19$ & $773$ & $15.0$ & $15$\\ \hline
H.B. Lawson & $2576$ & $27.4$ & $22$ & $1846$ & $23.2$ & $21$\\ \hline
P. Lax & $4601$ & $36.6$ & $30$ & $3560$ & $32.2$ & $27$\\ \hline
E. Lieb & $5147$ & $38.7$ & $33$ & $4349$ & $35.6$ & $32$\\ \hline
T. Liggett & $1975$ & $24$ & $17$ & $984$ & $16.9$ & $16$\\ \hline
L. Lovasz & $5638$ & $40.5$ & $34$ & $4259$ & $35.2$ & $30$\\ \hline
G. Lusztig & $5786$ & $41.1$ & $40$ & $4945$ & $38.0$ & $38$\\ \hline
R. MacPherson & $2031$ & $24.3$ & $22$ & $1676$ & $22.1$ & $21$\\ \hline
G. Margulis & $2267$ & $25.7$ & $26$ & $1788$ & $22.8$ & $25$\\ \hline
J. Mather & $1399$ & $20.2$ & $21$ & $1399$ & $20.2$ & $21$\\ \hline
B. Mazur & $2842$ & $28.8$ & $26$ & $2440$ & $26.7$ & $24$\\ \hline
D. McDuff & $2289$ & $25.8$ & $24$ & $1417$ & $20.3$ & $23$\\ \hline
H. McKean & $2480$ & $26.9$ & $24$ & $1701$ & $22.3$ & $21$\\ \hline
C. McMullen & $1738$ & $22.5$ & $25$ & $1368$ & $20.0$ & $24$\\ \hline
J. Milnor & $7856$ & $47.9$ & $48$ & $4559$ & $36.5$ & $40$\\ \hline
J. Morgan & $1985$ & $24.1$ & $25$ & $1484$ & $20.8$ & $21$\\ \hline
G. Mostow & $1180$ & $18.5$ & $19$ & $896$ & $16.2$ & $17$\\ \hline
J. Nash & $1337$ & $19$ & $10$ & $1337$ & $19$ & $10$\\ \hline
E. Nelson & $1010$ & $17.2$ & $15$ & $753$ & $14.8$ & $11$\\ \hline
L. Nirenberg & $9145$ & $51.6$ & $45$ & $8781$ & $50.6$ & $43$\\ \hline
\end{tabular}
\end{table}

\begin{table}[h]
\centering
\begin{tabular}{|c|c|c|c|c|c|c|}
\hline
Member & $N_{\tt citations}$ & Rule of thumb est. & $h$ & non-books only & revised est. & revised $h$\\
\hline
\hline
S. Novikov & $2368$ & $26.3$ & $27$ & $1677$ & $22.1$ & $21$\\ \hline
A. Okounkov & $1677$ & $24$ & $22.1$ & $1677$ & $24$ & $22.1$\\ \hline
D. Ornstein & $1100$ & $17.9$ & $19$ & $1022$ & $17.3$ & $18$\\ \hline
J. Palis & $1570$ & $21.4$ & $19$ & $895$ & $16.2$ & $18$\\ \hline
P. Rabinowitz & $6633$ & $44$ & $29$ & $5316$ & $39.4$ & $29$\\ \hline
M. Ratner & $506$ & $12.1$ & $13$ & $506$ & $12.1$ & $13$ \\ \hline
K. Ribet & $1022$ & $17.3$ & $18$ & $1021$ & $17.3$ & $18$\\ \hline
P. Sarnak & $3114$ & $30.1$ & $32$ & $2780$ & $28.5$ & $29$\\ \hline
M. Sato & $738$ & $14.7$ & $12$ & $738$ & $14.7$ & $12$\\ \hline
R. Schoen & $3945$ & $33.9$ & $34$ & $3493$ & $31.9$ & $34$\\ \hline
J. Serre & $10119$ & $54.3$ & $53$ & $4481$ & $36.1$ & $36$\\ \hline
%J. Serrin & $4539$ & $36.4$ & $37$ & $4430$ & $35.9$ & $36$\\ \hline
C. Seshadri & $984$ & $16.9$ & $15$ & $831$ & $15.6$ & $14$\\ \hline
Y. Sinai & $3357$ & $31.3$ & $31$ & $2547$ & $27.3$ & $28$\\ \hline
I. Singer & $2982$ & $29.5$ & $28$ & $2951$ & $29.3$ & $28$\\ \hline
Y. Siu & $1494$ & $20.9$ & $22$ & $1350$ & $19.8$ & $21$\\ \hline
S. Smale & $4581$ & $36.5$ & $39$ & $3942$ & $33.9$ & $36$\\ \hline
R. Solovay & $781$ & $15.1$ & $14$ & $781$ & $15.1$ & $14$\\ \hline
J. Spencer & $758$ & $14.9$ & $15$ & $1334$ & $19.7$ & $18$\\ \hline
R. Stanley & $6510$ & $43.6$ & $35$ & $3148$ & $30.3$ & $32$\\ \hline
H. Stark & $678$ & $14.1$ & $13$ & $653$ & $13.8$ & $13$\\ \hline
C. Stein & $763$ & $14.9$ & $12$ & $658$ & $13.9$ & $11$\\ \hline
E. Stein & $14049$ & $64$ & $49$ & $5788$ & $41.1$ & $37$\\ \hline
R. Steinberg & $1850$ & $23.2$ & $19$ & $1068$ & $17.6$ & $17$\\ \hline
S. Sternberg & $2438$ & $26.7$ & $25$ & $1476$ & $20.8$ & $18$\\ \hline
D. Stroock & $3299$ & $31.0$ & $27$ & $2028$ & $24.3$ & $24$\\ \hline
D. Sullivan & $3307$ & $31.1$ & $32$ & $3248$ & $30.8$ & $31$\\ \hline
R. Swan & $1109$ & $18$ & $20$ & $998$ & $17.1$ & $19$ \\ \hline
E. Szemer\'edi & $2536$ & $27.2$ & $26$ & $2536$ & $27.2$ & $26$\\ \hline
T. Tao & $6730$ & $44.3$ & $40$ & $6214$ & $42.3$ & $39$\\ \hline
J. Tate & $2979$ & $29.5$ & $30$ & $2612$ & $27.6$ & $28$\\ \hline
C. Taubes & $1866$ & $23.2$ & $24$ & $1626$ & $21.8$ & $23$\\ \hline
J. Thompson & $789$ & $15.2$ & $14$ & $789$ & $15.2$ & $14$\\ \hline
J. Tits & $3463$ & $31.8$ & $28$ & $2958$ & $29.4$ & $26$\\ \hline
K. Uhlenbeck & $1852$ & $23.2$ & $21$ & $1756$ & $22.6$ & $20$\\ \hline
S. Varadhan & $2894$ &  $29$ & $28$ & $2153$ & $25.1$ & $26$\\ \hline
D. Voiculescu & $2952$ & $29.3$ & $29$ & $2387$ & $26.4$ & $26$\\ \hline
A. Wiles & $1387$ & $20$ &  $14$ & $1387$ & $20$ &  $14$\\ \hline
S-T. Yau & $7536$ & $46.9$ & $44$ & $7066$ & $45.4$ & $43$\\ \hline
E. Zelmanov & $1055$ & $17.5$ & $16$ & $1020$ & $17.2$ & $16$\\ \hline
\end{tabular}

\end{table}

\end{document}